\documentclass[moor]{informs1hack}              % for a regular run
\makeatletter
\let\bfseries=\undefined
\DeclareRobustCommand\bfseries
{\not@math@alphabet\bfseries\mathbf
  \boldmath\fontseries\bfdefault\selectfont}
\makeatother

\newcommand{\Z}{{\mathbb Z}}

\newcommand{\R}{\mathbb R}

\def\Graver{{\mathcal G}}

\def\Circuits{{\mathcal C}}

\usepackage{enumerate}
\def\ve#1{\mathchoice{\mbox{\boldmath$\displaystyle\bf#1$}}
{\mbox{\boldmath$\textstyle\bf#1$}}
{\mbox{\boldmath$\scriptstyle\bf#1$}}
{\mbox{\boldmath$\scriptscriptstyle\bf#1$}}}
\newcommand\vea{{\ve a}}
\newcommand\veb{{\ve b}}
\newcommand\vecc{{\ve c}}

\newcommand\veg{{\ve g}}
\newcommand\veh{{\ve h}}
\newcommand\vel{{\ve l}}

\newcommand\ves{{\ve s}}
\newcommand\vet{{\ve t}}
\newcommand\veu{{\ve u}}
\newcommand\vev{{\ve v}}

\newcommand\vex{{\ve x}}

\newcommand\vez{{\ve z}}

\newcommand{\maO}{O}            % zero matrix

\usepackage{algorithm}

\usepackage{ifthen}
\newcommand{\DeclareBracket}[3]{
  \newcommand{#1}[2][]{%
  \ifthenelse%
  {\equal{##1}{}}%
  {\left#2##2\right#3}%
  {\csname ##1l\endcsname#2##2\csname ##1r\endcsname#3}}}
\DeclareBracket\set\{\}
\DeclareBracket\floor\lfloor\rfloor
\DeclareBracket\ceil\lceil\rceil
\DeclareBracket\bracket[]
\DeclareBracket\paren()

\newenvironment{psmallmatrix}{\left(\smallmatrix}{\endsmallmatrix\right)}

\newcommand\FourBlockBig[5][\relax]{\begin{pmatrix}#2& #3\\#4&#5 \end{pmatrix}\ifx#1\relax\else^{(#1)}\fi}
\newcommand\FourBlock[5][\relax]{\begin{psmallmatrix}#2& #3\\#4&#5 \end{psmallmatrix}\ifx#1\relax\else{^{(#1)}}\fi}
\newcommand\TwoBlock[3][\relax]{[#2,#3]\ifx#1\relax\else{^{(#1)}}\fi}

\newcommand{\red}{\sqsubseteq}
\newcommand{\T}{{\intercal}} %transpose

\usepackage{natbib}
 \NatBibNumeric
 \def\bibfont{\small}%
 \bibpunct[, ]{[}{]}{,}{n}{}{,}%

\usepackage[colorlinks=true,breaklinks=true,bookmarks=true,urlcolor=blue,
     citecolor=blue,linkcolor=blue,bookmarksopen=false,draft=false]{hyperref}

\def\EMAIL#1{\href{mailto:#1}{#1}}% When hyperref is used, otherwise outcomment
\def\URL#1{\href{#1}{#1}}         % When hyperref is used, otherwise outcomment

\TheoremsNumberedThrough     % Preferred (Theorem 1, Lemma 1, Theorem 2)
\EquationsNumberedThrough    % Default: (1), (2), ...
\begin{document}
%%%%%%%%%%%%%%%%

%\setcounter{page}{1}
%\DOI{10.1287/xxxx.0000.0000}%

 \RUNAUTHOR{De Loera, Hemmecke and Lee}

\RUNTITLE{Augmentation in Linear and Integer Linear Programming}

% \TITLE{Bundling Information Goods of Decreasing Value}
\TITLE{On Augmentation Algorithms for Linear \\ and Integer-Linear Programming:  \\ From Edmonds-Karp to Bland and Beyond}

\ARTICLEAUTHORS{%
\AUTHOR{Jes\'us A. De Loera}
\AFF{Department of Mathematics \\ University of California \\ One Shields Avenue \\ Davis, CA 95616 USA\\ \EMAIL{deloera@math.ucdavis.edu} \\ \URL{http://www.math.ucdavis.edu/$\sim$deloera}}
\AUTHOR{Raymond Hemmecke}
\AFF{Zentrum Mathematik \\ Technische Universit\"at M\"unchen \\ 85747 Garching bei M\"unchen, Germany \\ \EMAIL{hemmecke@tum.de} \\ \URL{http://www-m9.ma.tum.de/Allgemeines/RaymondHemmecke}}
\AUTHOR{Jon Lee}
\AFF{Department of Industrial and Operations Engineering \\ University of Michigan \\ 1205 Beal Avenue \\ Ann Arbor, MI 48109-2117 USA \\ \EMAIL{jonxlee@umich.edu} \\ \URL{https://sites.google.com/site/jonleewebpage/}}
} % end of the block

\ABSTRACT{
Motivated by Bland's linear-programming generalization of the renowned Edmonds-Karp efficient refinement of the Ford-Fulkerson maximum-flow algorithm, we discuss three closely-related natural augmentation rules for linear and integer-linear optimization. In several nice situations, we show that polynomially-many augmentation steps suffice to reach an  optimum. In particular, when using ``discrete steepest-descent augmentations'' (i.e., directions with the best ratio of cost improvement per unit 1-norm length), we show that the
number of augmentation steps is bounded by the number of elements in the Graver basis of the problem matrix,
giving the first ever strongly polynomial-time algorithm for $N$-fold integer-linear optimization. Our results also
improve on what is known for such algorithms in the context of linear optimization (e.g., generalizing the bounds of Kitahara and Mizuno for the number of
steps in the simplex method) and are closely related to research on the diameters of polytopes and the search for a strongly polynomial-time simplex
or augmentation algorithm.

%For the simplex method one is interested on the number of pivots needed to reach
%the optimum. In recent years several researchers have extended the notion of pivoting algorithm: For linear programs
%one can allow instead of just edges the entire range of circuits of the defining matrix.  The fundamental difference is
%that unlike the simplex pivoting we now allow moves that go through the interior of the polyhedral feasible region.
%Similarly, it has been proposed how to solve integer programs via similar "pivoting" or augmentation techniques,
%where the augmentation or pivoting steps are integer vectors from
%the Graver basis of the defining matrix.
%
%
%This paper contributes several bounds on the number of augmentation steps for
%three different pivot rules, and in both the context of LP and IPs.
%In some nice situations we show polynomially many augmentation steps suffice to reach the optimum. For instance, when using
%augmentation along directions with best ratio of cost improvement/unit length, we show that the
%number of augmentation steps is bounded by the number of elements in the Graver basis of the problem matrix, giving new strongly
%polynomial-time algorithms  of classical cases of strongly polynomial LPs (e.g., totally unimodular and for
%recovering the Edmonds-Karp algorithm for max flow ), but most important, giving strongly-polynomial time algorithms for some
%families  ILPs.
}

%\KEYWORDS{deterministic inventory theory; infinite linear programming duality;
%  existence of optimal policies; semi-Markov decision process; cyclic schedule}
%\MSCCLASS{Primary: 90B05; secondary: 90C40, 90C90}
%\ORMSCLASS{Primary: Inventory/production: deterministic multi-item;
%  secondary: dynamic programming/optimal control: deterministic
%  semi-Markov; programming: infinite dimensional}
%\HISTORY{Received November 20, 2003; revised March 8, 2004, and March 26, 2004.}

\KEYWORDS{augmentation, Graver basis, test set, circuit, elementary vector, linear program, integer program, Edmonds-Karp, steepest descent.}
\MSCCLASS{%
Primary: 90C10; %Integer programming
secondary:
65K05, %Mathematical programming methods
90C05, %Linear programming
52B55  %Computational aspects related to convexity
}%
\ORMSCLASS{Primary: Integer programming (theory);
secondary: Linear programming (theory); Mathematics (combinatorics, convexity).}
\HISTORY{}

\maketitle
%%%%%%%%%%%%%%%%%%%%%%%%%%%%%%%%%%%%%%%%%%%%%%%%%%%%%%%%%%%%%%%%%%%%%%

% Samples of sectioning (and labeling) in MOOR.
% NOTE: (1) all section levels end with a period,
%       (2) capitalization is as shown (sentence style, not title style).
%
%\section{Introduction.}\label{intro} %%1.
%\subsection{Duality and the classical EOQ problem.}\label{class-EOQ} %% 1.1.
%\subsection{Outline.}\label{outline1} %% 1.2.
%\subsubsection{Cyclic schedules for the general deterministic SMDP.}
%  \label{cyclic-schedules} %% 1.2.1
%\section{Problem description.}\label{problemdescription} %% 2.

\section{Introduction.}
We consider a general framework for solving linear programs (LPs) and integer-linear programs (ILPs) of the form
\begin{equation} \label{themilp}
  \min\set{\,\vecc^\T\vex\ :\ A\vex=\veb,\ \ve 0\leq\vex\leq\veu,\ \vex\in X\,},
\end{equation}
where $A\in\Z^{d\times n}$, $\veb\in\Z^d$, $\vecc\in\Z^n$, and where $X=\R^n$ (LP) or $X=\Z^n$ (ILP). We focus on solution algorithms that are based on an \emph{augmentation procedure}. At each iteration, we have a current feasible point $\vex_k$. An \emph{augmentation direction} $\vez$ has $\vex_k + \alpha\vez$ feasible for some $\alpha>0$ (which implies that
$\vez$ is in the kernel of $A$) and has $\vecc^\T \vez < 0$.
Together, the \emph{augmentation} is $\alpha\vez$, and we pass to the next
feasible solution  $\vex_{k+1}:=\vex_k + \alpha\vez$.
We note the trivial fact that $\vecc^\T \vex_{k+1} < \vecc^\T \vex_{k}$.
In this way, an augmentation procedure produces
a sequence of feasible solutions $\vex_1$, $\vex_2$, $\ldots$
that are successive improvements on the value of the objective function,
until an optimal solution is reached. Finally, an augmentation $\alpha\vez$ is \emph{maximal} if,
considering $\vex_k$ and $\vez$ to be fixed,
$\alpha$ is the largest value of $t>0$ for which  $\vex_{k+1}:=\vex_k + t\vez$ is feasible.
If the augmentation is maximal, $\vex_{k+1}$ is the best
point on the intersection of the feasible region with
the half-line $\{ \vex_k + t\vez ~:~ t >0\}$. In what follows, we give upper bounds on the number of maximal augmentations required to reached an optimum, for
both the LP and ILP cases.

Of course, augmentation algorithms are nothing new, and one of the best-known in the family is the classical  simplex algorithm. 
For the case of LP, the simplex method is indeed an augmentation algorithm, where
we start at a vertex of the polyhedron that is the feasible region,
and the augmentation used at each iteration corresponds to an available edge direction at the current vertex.
Of course, in the degenerate case, there can be considerable work to calculate an improving edge direction.
Nonetheless, by limiting augmentation directions to available edge directions at the current vertex, and always choosing maximal
augmentations, we insure that the next feasible solution is also a vertex, and so the simplex algorithm can continue.
Similarly, the idea of augmentation has played a very important historical role in the algorithmic theory of network flows; in particular, there is the seminal and very well-known work of  Edmonds and Karp (see \cite{edmondskarp}). They showed that for maximum-flow (essentially the same problem as maximizing the flow on a single arc of a flow-conservative network, subject to simple bounds on the other arcs), the number of augmentations (taken by the classic Ford-Fulkerson augmentation algorithm)
is bounded by the number of arcs times the number of vertices, or slightly more crudely by the number of arcs squared,
\emph{when augmentations are always chosen to have the fewest number of arcs, and the augmentation is maximal.}
In unpublished work, R. Bland (\cite{blandFlorida}) extended this, rather elegantly, to general LPs. Besides not being published
by Bland in the 1970's, it is not well known even today. The result was implicitly alluded to in print in 1987 (see \cite{BixMarTr}), and mentioned more concretely in J. Lee's 1986 dissertation \cite{LeeCornell} (see Proposition 3.1 in the follow-on publication \cite{LeeLAA}). Bland himself made a concrete statement of it (still without proof) in 1992  (see \cite{blandjensen}):

%%%%%%%%%%%%DONT TOUCH QUOTATION%%%%%%%%%%%%%%%%%%%%%%%%%%
\begin{quotation}
\it
``It was prompted by another result in the same Edmonds-Karp paper \cite{edmondskarp}, that if one always augments on a shortest augmenting
path, the number of augmentations in the Ford-Fulkerson maximum flow algorithm
is less than the product of the numbers of nodes and edges. Fulkerson \cite{fulk}  had
investigated the extent to which fundamental properties of networks generalize to
broader classes of linear programming problems, where elementary vectors in the
appropriate subspaces play the roles of circuits and cocircuits (minimal cutsets).
Bland's dissertation \cite{blandthesis} carried this further, and in later work he showed how the
Edmonds-Karp result generalizes to arbitrary linear programming problems of the
form
{\rm
\begin{eqnarray}
&&\mbox{maximize }  x_0 \nonumber \\
&&\mbox{subject to } A\vex =\ve 0, \label{blandlp}\\
&&\vel\leq \vex \leq \veu. \nonumber
\end{eqnarray}
}
Here the bound on the number of augmentations is the number of variables times
the number of different lengths of normalized elementary vectors; in particular if
$A$ is totally unimodular the bound is the product of the dimensions of $A$, as in the
max flow result of Edmonds and Karp, where $A$ is the (totally unimodular) node-edge
matrix of a directed graph. This seemed to be amusing, but without any obvious
use, until late in 1978 when Paul Seymour \cite{seymour} proved his remarkable decomposition
theorem for unimodular matroids. Bland and Edmonds \cite{blandedmonds} used Seymour's theorem
to show how shortest augmentations could be computed in polynomial time when
the constraint matrix $A$ is totally unimodular. Cost scaling was used to extend the
approach to a polynomial-time primal algorithm for totally unimodular linear programming problems
{\rm
\begin{eqnarray*}
&&\mbox{maximize }  \vea^\T\vex  \\
&&\mbox{subject to } A\vex =\veb, \\
&&\vel\leq \vex \leq \veu.
\end{eqnarray*}
}
which are solved as a sequence of subproblems of the form (\ref{blandlp}). More recently
Tardos \cite{Tardos} has given strongly polynomial algorithms for totally unimodular linear
programming problems, and still more general classes of combinatorial optimization
problems, using cost scaling in a far more clever way.''
\end{quotation}
%%%%%%%%%%%%%%%%%%%%%%%%%%%%%%%%%%%%%%%%%%

Today, progress on the augmentation approach to LP continues: new families of augmentation algorithms for
LP include those proposed in \cite{cardosoclimaco} (and rediscovered in \cite{baraszvempala}). Very recently, 
the work of Kern, Faigle, and Peis \cite{peis} rekindled our interest in Bland's result as they, among other things, rediscovered it. Additionally,  Kitahara and Mizuno  \cite{kitaharamizuno} gave an upper bound for  the number of different basic-feasible solutions generated by the simplex method. Their results are direct generalizations of Y. Ye's work \cite{yye}.

%Former students of Bland have seen and at least one (the third author of the present paper) has even taught the nice proof of Bland's
%result, but unfortunately it never appeared anywhere.

On the other hand, the study of augmentation algorithms for ILPs goes back at least to the 1970's and was part of the study of 
``test sets'' (see \cite{Graver:75,Scarf:97, SchulzWeismantel99} and the references therein).  More recently 
there have been important developments. In \cite{SW2002} the authors show that one can solve every integer-linear programming problem in polynomial time provided one can efficiently solve a  special \emph{directed-augmentation problem} (the directed-augmentation problem differs from the ordinary augmentation problem by splitting each direction into its positive and negative parts and considering linear objectives on each of these parts). These authors showed that if one can solve the directed-augmentation problem in polynomial time, then the original problem can be solved by a polynomial-time algorithm. Their main application is to specific combinatorial-optimization problems, such as the min-cost flow problem. Later \cite{Hemmecke+Onn+Weismantel:oracle} showed that if a ``best-improving'' augmentation $\vez$ in the ``Graver basis'' of the constraint matrix is chosen (the authors of \cite{Hemmecke+Onn+Weismantel:oracle} called this augmentation rule  ``greedy improvement''; here we will call this ``deepest descent''), only polynomially-many augmentations (in the binary encoding length of the input data) must be performed in order to reach the optimal value. The paper \cite{Hemmecke+Onn+Weismantel:oracle} shows how to efficiently use an augmentation algorithm for separable convex objective functions and implicitly provided bounds like those presented here.

%%%%%%%%%%%%%%%%%%%%%%%%%%%%%%%

The set of possible augmentation directions we use depends on the type of problem:
for the LP case, we use as the set of possible augmentation directions the \emph{circuits} of $A$. $\Circuits(A)$ consists of the normalized \emph{elementary vectors} (or \emph{circuits}) associated with $\ker(A)\setminus\set{\,\ve 0\,}$ (see \cite{rocka}) --- that is, the vectors having (set-wise)  minimal support in $\ker(A)\setminus\set{\,\ve 0\,}$. The set of elementary vectors of $\ker(A)\setminus\set{\,\ve 0\,}$ is a finite set of lines through the origin, with the origin excluded. Usually, it is convenient to normalize in an arbitrary manner, so that $\Circuits(A)$ comprises a single point and its negative from each such line. In our context, Bland's normalization uses the objective function so that $\vecc^\T \vez=-1$ for every augmentation direction --- in his terms, such a $\vez$ is \emph{unit augmenting}. For the vectors in $\Circuits(A)$, we choose on each line the (nonzero) integer point closest to the origin and its negative as the normalized representatives.

For the IP case, the set of possible augmentation directions is  the \emph{Graver basis} of $A$, denoted by $\Graver(A)$. We obtain Graver's original finite set  of $\red$-minimal elements in $\ker(A)\cap\Z^n\setminus\set{\,\ve 0\,}$, where $\veu\red\vev$ if and only if $u^{(i)}v^{(i)}\geq 0$ and $|u^{(i)}|\leq|v^{(i)}|$ (see \cite{Graver:75}).
In general, $\Graver(A)$ has a nice sign-compatible representation property: every (integer vector) $\vez\in\ker(A)$ can be written as $\vez=\sum \alpha_i\veg_i$, with $\veg_i\in\Graver(A)$, $\alpha_i>0$, $\alpha_i\in\Z$, and $\alpha_i\veg_i\red\vez$, for all $i$, and the sum involves at most $2n-2$ terms.
 In fact, $\Graver(A)$ is an inclusion-minimal set with this property. It should be noted that due to the sign-compatible representation we have
\[
  \|\vez\|_1=\sum\alpha_i\|\veg_i\|_1,
\]
where $\|\cdot\|_1$ is the usual 1-norm.

We remark that the all circuits are members of the Graver basis. But they also provide an elegant sign-compatible representation 
\emph{over the reals} for all non-zero elements of $\ker(A)$. Specifically, every vector $\vez\in\ker(A)$ can be written as 
$\vez=\sum \alpha_i\veg_i$, with $\veg_i\in{\cal C}(A)$, $\alpha_i\geq0$, $\alpha_i\in\R$, and the sum involves at most $n$ terms. Note that the
(worst-case) number of summands is smaller than in the general Graver decomposition. It is worth remarking that Graver bases can also be defined 
for the general mixed-integer case (see \cite[Section 3.3]{Hemmecke:PSP}) --- but in the mixed-integer setting, the Graver basis is not generally finite, 
in sharp contrast to the LP and ILP cases; so we do not consider that situation here as, unlike the LP and ILP cases, it requires the existence of some 
kind of specialized oracle to generate the augmentations one by one. For an introduction to Graver bases and the latest on augmentation algorithms for integer
programming see the books \cite{DHKbook} and \cite{onn:nonlinear-discrete-monograph}.

%%%%%%%%%%%%%%%%%%%%%%%%%%%%%%%

\subsection{Our Contributions}

The main contribution of this paper is an extension of Bland's theorem beyond LP to ILP.
Our main results and proofs extend what Bland did for circuits to Graver bases. We look at several  augmentation rules.  
Although they are very similar, our results here are naturally divided into the LP and ILP cases, with small but important 
technical changes. Also, depending on whether variables are restricted to be integers or not, we rely on different sets as the set 
of allowable augmentation directions, which we denote by $\cal T (A)$. The set $\cal T(A)$ will change depending on the case.

In what follows, we will employ $\cal T(A)$ in augmentation algorithms that iteratively replace a feasible (either continuous or integer)
solution $\vex_k$ to $A\vex=\veb$, $\ve 0\leq\vex\leq\veu$, by a better feasible solution $\vex_{k+1}:=\vex_k+\alpha\vez$, where
 $\vecc^\T\vez<0$ and $\alpha>0$. We consider three specific augmentation rules.

\medskip
  \begin{definition}[Discrete Deepest Descent]
    With respect to a feasible solution $\vex_k$, we choose  $\vez$ such that $-\alpha\vecc^\T\vez$ is maximized among all $\vez\in{\cal T}(A)$ and $\alpha>0$ such that $\vex_{k+1} :=\vex_k + \alpha \vez$ is feasible (note that for ILP this means that $\alpha\in\Z$).
  \end{definition}

\medskip

  \begin{definition}[Discrete Dantzig Descent]
    With respect to a feasible solution $\vex_k$, we choose  $\vez$ such that $-\vecc^\T\vez$ is maximized among all $\vez\in\cal T(A)$ such that $\vex_k + \epsilon\vez$ is feasible for some $\epsilon>0$ (note that for LP this means for all sufficiently small $\epsilon>0$, and for ILP this means for $\epsilon=1$). Then we take a maximum augmentation in such a direction. That is, we let $\vex_{k+1} := \vex_k + \alpha \vez$, where $\alpha$ is the largest value for which $\vex_k + \alpha \vez$ is feasible (note that for ILP this means that $\alpha\in\Z$).
  \end{definition}

\medskip

  \begin{definition}[Discrete Steepest Descent]
    With respect to a feasible solution $\vex_k$, we choose  $\vez$ such that $-\vecc^\T\vez/\|\vez\|_1$  is maximized among all $\vez\in\cal T(A)$ such that $\vex_k + \epsilon \vez$ is feasible for some $\epsilon>0$ (note that for LP this means for all sufficiently small $\epsilon>0$, and for ILP this means for $\epsilon=1$). Then we take a maximum augmentation  in such a direction. That is, we let $\vex_{k+1} := \vex_k + \alpha \vez$, where $\alpha$ is the largest value for which $\vex_k + \alpha \vez$ is feasible  (note that for ILP this means that $\alpha\in\Z$).
  \end{definition}

\medskip

Note that all three augmentation rules produce maximal augmentations. Practically speaking, there are more situations in which discrete steepest descent and 
discrete Dantzig descent can be practically implemented, as compared to discrete deepest descent. Still, it is interesting to analyze and contrast these augmentation rules.
We derive the following main results concerning these augmentation rules for the ILP and the LP cases. In some structured situations, these bounds will provide very 
good guarantees of performance. Our ILP theorem is new, while our LP results extend what Bland started early on, and greatly extends the applicability of the 
bounds of Kitahara and Mizuno. 

Next we state the results. Note that in what follows the base of all logarithms is two.

For the discrete deepest-descent rule, in the ILP case a polynomial time bound was first proved, but not explicitly stated as a proposition, in \cite{Hemmecke+Onn+Weismantel:oracle} and in \cite{onn:nonlinear-discrete-monograph} (see top of page 47 (end of proof of Lemma 3.10)).  Here we make it very explicit:

\begin{lemma}[(see also \cite{Hemmecke+Onn+Weismantel:oracle})] \label{ILP-greedy}
The number of discrete deepest-descent augmentations needed to reach an optimal solution from $\vex_0$ is bounded by $(4n-4)\,\log(\vecc^\T (\vex_0-\vex_{\min})).$
\end{lemma}

In addition, for the other two augmentation rules we prove:

\begin{theorem}[ILP case]\label{theorem:main, ILP case}
Let $A\in\Z^{d\times n}$, $\veb\in\Z^d$ and $\vecc\in\Z^n$ define the ILP
\[
\min\set{\,\vecc^\T\vex\ :\ A\vex=\veb,\ \ve 0\leq\vex\leq\veu,\ \vex\in\Z^n\,}.
\]
Let $\vex_0$ be an initial feasible solution, let $\vex_{\min}$ be an optimal solution, and let $\gamma$ be the maximum non-zero 
entry (in absolute value) in any feasible solution. Then we have the following bounds on the number of augmentations to reach an optimal solution from $\vex_0$.
\begin{enumerate}
\medskip
\item[(a)] The number of discrete Dantzig-descent augmentations needed to reach an optimal solution of the ILP is bounded by 
$(4n-4)\gamma \, \log(\vecc^\T (\vex_0-\vex_{\min})).$
%$O(n\,\gamma\,\log(\vecc^\T (\vex_0-\vex_{\min})))$.
\item[(b)] Any discrete steepest-descent direction (which by definition belongs to $\Graver(A)$) is an overall steepest-descent direction (which could be any applicable direction from $\Z^n$). Moreover, the number of discrete steepest-descent augmentations to reach an optimal solution of the given ILP is bounded by $|\Graver(A)|$.
\end{enumerate}
\end{theorem}

Next we present our results for augmentation algorithms in linear programming.
The results for \emph{integer} augmentations arguments in \cite{Hemmecke+Onn+Weismantel:oracle} and in \cite{onn:nonlinear-discrete-monograph} (see top of page 47 (end of proof of Lemma 3.10)) could be adapted to obtain an LP bound too, but the bound one obtains from  \cite{Hemmecke+Onn+Weismantel:oracle}  leads directly to an extra factor of $n$ in the number of augmentations due to the strategy presented there to ensure termination. Our contribution is that we manage to get rid of this factor of $n$:

\begin{lemma} \label{LP-greedy}
 The number of discrete deepest-descent augmentations needed to reach an optimal solution from $\vex_0$ is bounded by $2n\log(\delta\,\vecc^\T (\vex_0-\vex_{\min})).$ 
\end{lemma}

Moreover we prove the following extensions to the other augmentation rules for LPs:

\begin{theorem}[LP case]\label{theorem:main, LP case}
Let $A\in\Z^{d\times n}$, $\veb\in\Z^d$ and $\vecc\in\Z^n$ define the LP
\[
\min\set{\,\vecc^\T\vex\ :\ A\vex=\veb,\ \ve 0\leq\vex\leq\veu,\ \vex\in\R^n\,}.
\]
Let $\vex_0$ be an initial feasible solution, let $\vex_{\min}$ be an optimal solution, let $\gamma$ be the maximum non-zero entry (in absolute value) in any feasible solution, and let $\delta$ denote the least common multiple of all subdeterminants of $A$. Then we have the following bounds on the number of augmentations to 
reach an optimal solution from $\vex_0$.
\begin{enumerate}
\medskip
\item[(a)] The number of discrete Dantzig-descent augmentations needed to reach an optimal solution from $\vex_0$ is no more than
$2n^2\delta\gamma\,\log(\delta\,\vecc^\T (\vex_0-\vex_{\min})).$
%$O(n^2\delta\gamma\,\log(\delta\,\vecc^\T (\vex_0-\vex_{\min})))$.
\item[(b)] Any discrete steepest-descent direction (which by definition belongs to ${\cal C}(A)$) is an overall steepest-descent direction (which could be any applicable direction from $\R^n$). Moreover, the number of discrete steepest-descent augmentations to reach an optimal solution of the given LP is bounded by $|{\cal C}(A)|$.
\end{enumerate}
\end{theorem}

While by nature LP augmentation algorithms cannot have cycling, as can happen for the simplex method, it should be noted that the LP case of the discrete deepest-descent and the discrete Dantzig-descent augmentation algorithms only guarantee that we reach a feasible solution $\vex_k$ with objective value ``close enough'' to the optimal value. Whether any given feasible solution is ``close enough'' can be checked by generating a vertex of the feasible region with objective value that is at least as good as $\vecc^\T\vex_k$. Finding such a vertex can be done by augmenting $\vex_k$ only along such (discrete deepest-/Dantzig-descent) circuit directions that lead to a better feasible solution with an additional component reaching its lower or upper bound. Geometrically, this corresponds to (iterative) augmentation within the smallest face of the polyhedron $\set{\,\vex\ :\ A\vex=\veb,\ \ve 0\leq\vex\leq\veu,\ \vex\in\R^n\,}$ that contains $\vex_k$. (The circuits of the problem matrix $A$ provide such a restricted optimality certificate). We recommend  the example in \cite{Hemmecke:PSP}. Finally, again using the circuits of $A$, it can be checked whether the vertex found is optimal. The overall closeness-test requires at most $n$ augmentations.

Finally we present the interesting consequences of the main theorems.
Note that the bounds in part (c) of Theorems \ref{theorem:main, ILP case} and \ref{theorem:main, LP case} only depend on $A$, but not on $\veb$, $\vecc$, and the particular initial feasible solution $\vex_0$ chosen. As a direct consequence of Part (b) of Theorem \ref{theorem:main, ILP case}, we obtain the following corollaries.

\begin{corollary} \label{01boxedcase}
  For the pure 0/1 ILP $\min\set{\,\vecc^\T\vex\ : A\vex=\veb,\ {\ve 0}\leq\vex\leq {\ve 1},\ \vex\in\Z^n\,}$ the number of discrete deepest-/Dantzig-/steepest-descent augmentations is bounded by $O(n \log(\|\vecc \|_1))$.
\end{corollary}

We note that in the LP case, the proof of Theorem \ref{theorem:main, LP case}, part (c), immediately gives Bland's
fundamental result:

\begin{corollary}[Bland's Theorem] \label{Corollary: Blands result}
  The number of discrete steepest-descent augmentations needed to solve   $\min\set{\,\vecc^\T\vex:A\vex=\veb, \ve 0\leq\vex\leq\veu, \vex\in\R^n\,}$ is bounded by the number of (different) positive values of $-\vecc^\T \vez/\|\vez\|_1$ over all (elementary vectors) $\vez\in{\cal C}(A)$ times the number $n$ of variables.
\end{corollary}

As for totally-unimodular matrices $A$, ${\cal C}(A)=\Graver(A)$, i.e., they coincide for the LP and ILP cases, the proof of Theorem \ref{theorem:main, LP case}, part (c), also implies:

\begin{corollary}\label{Corollary: TU matrix implies Blands result also for ILP}
  For totally-unimodular matrix $A$, the number of discrete steepest-descent augmentations needed to solve $\min\set{\,\vecc^\T\vex:A\vex=\veb, \ve 0\leq\vex\leq\veu, \vex\in\Z^n\,}$ is bounded by the number of (different) positive values of $-\vecc^\T \vez/\|\vez\|_1$ over all (elementary vectors) $\vez\in{\cal C}(A)$ times the number $n$ of variables.
\end{corollary}

For a totally-unimodular matrix $A$, $\Circuits(A)$ consists only of vectors with (at most $d+1$) entries in $\set{\,-1,0,1\,}$ Thus, for $\vez\in\Circuits(A)$ we have $|-\vecc^\T \vez|\leq\|\vecc\|_1$ and $\|\vez\|_1$ takes on at most $d+1$ different values. Plugging this into Corollary \ref{Corollary: TU matrix implies Blands result also for ILP} we get the following.

\begin{corollary} \label{tunumsteps}
  For totally-unimodular matrix $A$, the number of discrete steepest-descent augmentations needed to solve $\min\set{\,\vecc^\T\vex:A\vex=\veb, \ve 0\leq\vex\leq\veu, \vex\in\Z^n\,}$ is bounded by $n(d+1)\|\vecc\|_1$.
\end{corollary}

From this we immediately recover the complexity bound for the algorithm of Edmonds and Karp to find maximum flows in a network: let $A$ be the node/arc-incidence matrix of a connected directed graph. Note that one row of $A$ is linearly dependent on the other rows and thus can be removed from $A$. Then $n=|E|$ and $d=|V|-1$. As we maximize the flow on a specific (auxiliary) arc (from sink to source), we have $\|\vecc\|_1=1$. Thus, $n(d+1)\|\vecc\|_1=|E|\cdot|V|$ bounds the number of discrete steepest-descent augmentations to solve the max-flow problem. This observation is not  surprising, as the augmentation approach using Graver bases specializes to the algorithm by Edmonds and Karp in the setting of maximum flows.

%%%%%%%%%%%%%%%
%${\bf -----------------PLEASE \quad READ------------}$

It is  worth recalling that although the complexity statements in Corollary \ref{01boxedcase} and Corollary \ref{tunumsteps}
depend on the unary size of $\vecc$, these two results actually can be improved based on the results of Frank and Tardos \cite{andras+eva}.
Frank and Tardos used  Diophantine approximation to replace $\vecc$ with a new objective function $\vecc'$ where the integer
numbers occurring in the entries are small but define an equivalent problem with the same optima. Furthermore, the
new weight function $\vecc'$ can be found in strongly polynomial-time. In conclusion, if we are able to generate in polynomial time
the corresponding augmentation elements of the Graver basis according to one of the three rules, we  obtain strongly polynomial-time
algorithms. Of course, this is in general hard to do other than by computing the entire Graver basis, but for special matrices $A$, one can
do much better.

Our results on augmentations for totally-unimodular matrices provide another interesting geometric result. For years, researchers have been looking at the diameter
of the graph of polyhedra (i.e., the graph whose nodes are the vertices of the convex polyhedron $P=\{ \vex :  A\vex=\veb, \ve 0\leq\vex\leq\veu \}$). It has been shown in
\cite{determinants+diameters} that the diameter of the graph of totally-unimodular $d$-dimensional polytopes (i.e., the length of the longest shortest path between a pair of nodes) is bounded above by $d^{3.5}\log(d)$.  
We wish to stress that in our circuit augmentations, we do not always follow edges of the polyhedron $P$. Rather, we may cut through the interior of $P$ or the interior of some faces. This suggests the notion of \emph{circuit diameter}. The circuit distance from $v_1$ to $v_2$  is the smallest number of circuit augmentations needed to go from $v_1$ to $v_2$. We can then define the circuit diameter   as the maximum number of steps along circuit basis directions that are needed to go from any vertex of the polyhedron to any other vertex of the polyhedron. This notion was first introduced and investigated in the article \cite{circuitdiameter}. More recently other  generalizations of the notion of diameter 
using circuits were introduced in  \cite{manycircuitdiameters}. One can show that the circuit diameter of a polyhedron is bounded from above by the usual combinatorial diameter of polyhedra (this is because arcs are themselves  some of the possible augmentation directions, a subset of the circuits). In this context, our results show that one can bound the \emph{circuit diameter} for totally-unimodular polytopes in standard form as follows  %By taking $\vecc$ equal to a unit $0/1$ vector we obtain.

\begin{corollary}
  For a $d \times n$ totally-unimodular matrix $A$, the circuit diameter of the polyhedron $P:=\set{\,\vex :A\vex=\veb, \ve 0\leq\vex\leq\veu}$ is bounded above by $2(n(d+1)(n-d))$.
\end{corollary}

This corollary is a good general bound for totally-unimodular matrices, but we suspect it can be further improved, as we already know that for network-flow polytopes a better bound is possible. Indeed, Orlin \cite{orlin97-networkdiam} designed a polynomial-time primal network-simplex algorithm for the minimum-cost flow problem which gives a graph-diameter bound of $O(|E||V| \log(|V|))$, with $E,V$ equal to the sets of of arcs and nodes in the network. 

%----------------------------------
\proof{Proof:} To see this, take any vertex $v$ of $P$. Choose a cost vector specific to $v$: namely, let $c_i=1$ if $v_i=0$ and $c_i=-1$ if $v_i=u_i$. Otherwise, put $c_i=0$. Thus, this new objective function $\vecc$ has $n-d$ = $n-rank(A)$ many nonzero entries. Thus, by Corollary \ref{tunumsteps}, with the steepest-descent rule, any other vertex is connected to $v$ with no more than $n(d+1)(n-d)$ many augmentations. Hence,  the path between any pair of vertices of P is bounded by $2n(d+1)(n-d)$. \Halmos
\endproof
%%%%%%%%%%%%%%%%%%

Finally, let us consider optimization problems for which the constraint matrix is structured. Recall that an $N$-fold matrix is a matrix of the form
\[
  [A,B]^{(N)}:=\begin{pmatrix}
    B & B & \cdots & B \\
    A & \maO  &   & \maO  \\
    \maO  & A &   & \maO  \\
       &   & \ddots &   \\
    \maO  & \maO  &   & A
  \end{pmatrix}.
\]
For fixed matrices $A$ and $B$, the size of the Graver basis of $[A,B]^{(N)}$ (and its binary encoding length) increases only polynomially with $N$. Combining this with Theorems \ref{theorem:main, ILP case} and \ref{theorem:main, LP case}, parts (c), we obtain the following result.

\begin{theorem}\label{Theorem: N-fold LPs and IPs are strongly poly-time solvable}
  For fixed matrices $A$ and $B$, the associated families of $N$-fold LPs and ILPs can be solved in strongly polynomial time.
\end{theorem}

This generalizes the results from \cite{hemmecke-onn-romanchuk:nfold-cubic,Hemmecke+Onn+Weismantel:oracle}, which showed that for fixed matrices $A$ and $B$, the corresponding $N$-fold ILPs could be solved in time polynomial in $N$, in fact in $O(N^3)$ steps. Theorem \ref {Theorem: N-fold LPs and IPs are strongly poly-time solvable} strengthens this to strong polynomiality. As a direct consequence we obtain the following results:

\begin{corollary} The following special cases of $N$-fold matrices can be solved in strongly polynomial time in the linear and integer case:
  \begin{itemize}
    \item All $2$-way transportation problems with a fixed number of rows or columns.
    \item All $3$-way transportation problems with two of the three dimensions fixed.
    \item All $d$-way transportation problems with $d-1$ dimensions of fixed constant value.
  \end{itemize}
\end{corollary}

%%%%%%%%%%%%%%%%%%%%%%%%%%%%%%%%%%%%%%%%%%%%%%%%%
%%%%%%%%%%%%%%%%%%%%%PROOFS   %%%%%%%%%%%%%%%%%%%%%%%
%%%%%%%%%%%%%%%%%%%%%%%%%%%%%%%%%%%%%%%%%%%%%%%%%

\section{Proofs}

We present the proofs to the above results in separate subsections. We will use the following  lemma (a slight variation from Theorem 3.1 in \cite{Ahuja+Magnanti+Orlin}) which  establishes the bounds we claim once we can guarantee sufficient improvement at each iteration:
\begin{lemma}
\label{Lemma: Sufficient improvement leads to polynomial algorithm}

Let $\epsilon>0$ be given. Moreover, let $H$ denote the difference between maximum and minimum objective-function values of the LP/ILP problem in $n$ variables.
Suppose that $f^k=\vecc^\T \vex_k$ is the objective-function value of the solution $\vex_k$ at the $k$-th iteration of an algorithm and that $f^*=\vecc^\T \vex_{\min}$ is the minimum objective-function value. Furthermore, suppose that the algorithm guarantees that for every augmentation $k$,
\[
(f^k-f^{k+1})\geq\beta(f^k-f^*)
\]
(i.e., the improvement at augmentation $k+1$ is at least $\beta$ times the maximum possible improvement). Then the algorithm reaches a solution with $f^k-f^*<\epsilon$ in no more than $2 \log{(H/\epsilon)}/\beta$ augmentations.
\end{lemma}

\proof{Proof:} Without loss of generality, we assume that $\vecc$ is an integer vector. (Thus for the ILP version the values $f^k$ are decreasing successive integer values.) The biggest necessary change of value of the objective function occurs when we start the augmentation at an maximum point $\vex_{\max}$ and end with some minimum optimal point $\vex_{\min}$. Thus, we have $H=f^{\max}-f^{\min}$.

If we had at every augmentation an improvement of at least $\beta (f^{\max}-f^{\min})/2$, then
in no more than $2/\beta$ augmentations we would have reached the optimum. But if this improvement is not achieved at each augmentation, say at the $q$-th augmentation we have $f^q-f^{q+1} \leq \beta (f^{\max}-f^{\min})/2$ then, together with the hypothesis $(f^q-f^{q+1})\geq\beta(f^q-f^{\min})$, we get that
\[
  f^q-f^{\min} \leq (f^{\max}-f^{\min})/2=H/2.
\]
In other words, the overall improvement reached so far is at least half of the maximum possible improvement $H$.
In conclusion, after $2/\beta$ augmentations, we have either reached the optimum or have at least divided the possible gap by $2$. Therefore in no more than $2/\beta \log_2{(H/\epsilon)}$ augmentations we reach a solution with $f^k-f^*<\epsilon$. \Halmos
\endproof

It is important to note that in the ILP case $f^k$ is integer, and we can apply Lemma \ref{Lemma: Sufficient improvement leads to polynomial algorithm} with $\epsilon=1$ and conclude that we can reach the optimum in $O(\log{(H)}/\beta)$ augmentations. In the LP case, let $\delta$ denote the least common multiple of all subdeterminants of $A$. Observe that once we find a feasible solution $\vex_k$ with objective value $f^k$ satisfying $f^k-f^*<\epsilon=1/\delta$, then any vertex with an objective value of at most $f^k$ must be optimal. As explained above (right before the statement of Theorem \ref{theorem:main, ILP case}), such a vertex can be found from $\vex_k$ in at most $n$ additional augmentations. This leads to extra factors of $n$ and of $\log(\delta)$ in the bounds for the LP cases compared to the ILP cases.

\subsection{Proof of Lemma \ref{ILP-greedy} and Theorem \ref{theorem:main, ILP case}}

Let us assume that $\vex_k$ is a non-optimal feasible solution, and let $\vex_{\min}$ be an optimal solution to the ILP. Then there exists a (sign-compatible) representation
\[
  \vex_{\min} -  \vex_k =\sum\alpha_i\veg_i,
\]
with $\alpha_i>0$, $\alpha_i\in\Z$ and with $\alpha_i\veg_i\red \vex_{\min} - \vex_k$. Moreover, due to Seb\"o's result \cite{Sebo90}, at most $2n-2$ summands are needed.

Note that sign-compatibility of the representation $\vex_{\min}-\vex_k =\sum\alpha_i\veg_i$ implies that for all $i$ the vectors $\vex_{k}+\alpha_i\veg_i$ and $\vex_{\min}-\alpha_i\veg_i$ are all feasible solutions, since their components lie between the components of $\vex_k$ and of $\vex_{\min}$. Moreover, we can observe that for all such sign-compatible representations $\vex_{\min}-\vex_k =\sum\alpha_i\veg_i$ we must have $\vecc^\T\veg_i\leq 0$ for all $i$, as otherwise $\vex_{\min}-\alpha_i\veg_i$ would be a feasible solution with $\vecc^\T (\vex_{\min}-\alpha_i\veg_i)=\vecc^\T\vex_{\min}-\alpha_i\vecc^\T\veg_i<\vecc^\T\vex_{\min}$, contradicting the minimality of $\vex_{\min}$.

Next, we analyze what happens for each choice of augmentation rule:

\subsubsection{Proof of Lemma \ref{ILP-greedy}: Discrete deepest descent.} We observe that
\[
  0>\vecc^\T (\vex_{\min}-\vex_k)
   =\vecc^\T \sum\alpha_i\veg_i
   = \sum\alpha_i\vecc^\T\veg_i
   \geq -(2n-2)\Delta
\]
where $\Delta>0$ is the largest value of $-\alpha\vecc^\T\vez$ over all $\vez\in\Graver(A)$ and integer $\alpha>0$ for which $\vex_k+\alpha\vez$ is feasible. Rewriting this, we get
\[
  \Delta\geq\frac{\vecc^\T (\vex_k-\vex_{\min})}{2n-2}.
%\alpha_{\max}}.
\]
Now let $\alpha\vez$ be the discrete deepest-descent augmentation applied to $\vex_k$, leading to $\vex_{k+1}:=\vex_k+\alpha\vez$. Then we get $\Delta=-\alpha\vecc^\T\vez$ and
\[
  \vecc^\T (\vex_k-\vex_{k+1})=-\alpha\vecc^\T\vez=\Delta
  \geq\frac{\vecc^\T (\vex_k-\vex_{\min})}{2n-2}.
\]
Thus, we have a factor of $\beta=1/(2n-2)$ of objective-function decrease at each augmentation, leading to the desired polynomial number of augmentations via Lemma \ref{Lemma: Sufficient improvement leads to polynomial algorithm} taking $\epsilon=1$.
In this case, we get the number of augmentations bounded by $(4n-4)\,\,\log(\vecc^\T (\vex_0-\vex_{\min})))$.

\subsubsection{Proof of part (a) of Theorem \ref{theorem:main, ILP case}: Discrete Dantzig descent.} We observe that
\[
  0>\vecc^\T (\vex_{\min}-\vex_k)
   =\vecc^\T \sum\alpha_i\veg_i
   = \sum\alpha_i\vecc^\T\veg_i
   \geq - \Delta_0 \sum\alpha_i
   \geq -(2n-2)\Delta_0 \alpha_{\max},
\]
where this time $\Delta_0>0$ denotes the greatest value of $-\vecc^\T\vez$ over all $\vez\in\Graver(A)$ for which $\vex_k+\vez$ is still feasible and where $\alpha_{\max}=\max\set{\, \alpha_i\, }$. Rewriting this, we get
\[
  \Delta_0 \geq\frac{\vecc^\T (\vex_k-\vex_{\min})}{(2n-2) \alpha_{\max}}.
\]
Now let $\alpha\vez$ be the discrete Dantzig-descent augmentation  applied to $\vex_k$, leading to $\vex_{k+1}:=\vex_k+\alpha\vez$. Then we get
\[
  \vecc^\T (\vex_k-\vex_{k+1})=-\alpha\vecc^\T\vez=\alpha\Delta_0\geq\Delta_0
  \geq\frac{\vecc^\T (\vex_k-\vex_{\min})}{(2n-2) \alpha_{\max}}
  \geq \frac{\vecc^\T (\vex_0-\vex_{\min})}{(2n-2)\gamma},
\]
where $\gamma$ is the maximum entry in any feasible integer solution (or, equivalently, in any vertex of $P_I$). Thus, we have a factor of $\beta=1/((2n-2) \gamma)$ of objective-function decrease at each augmentation leading to the desired polynomial number of augmentations via Lemma \ref{Lemma: Sufficient improvement leads to polynomial algorithm} taking $\epsilon=1$. In this case we get the number of augmentations bounded by  $(4n-4)\,\gamma \, \log(\vecc^\T (\vex_0-\vex_{\min})))$

%\subsubsection{Proof to part (c): Discrete steepest descent.}  Here we have
%\[
%  0>\vecc^\T (\vex_{\min}-\vex_k)=\vecc^\T \sum\alpha_i\veg_i = \sum \alpha_i\|\veg_i\|_1
%  \left(\frac{\vecc^\T\veg_i}{\|\veg_i\|_1}\right)
%  \geq -\Delta \left\|\sum\alpha_i\veg_i\right\|_1
%   = -\Delta \|\vex_{\min}-\vex_k\|_1
%  \geq -2n\gamma\Delta,
%\]
%where $\Delta>0$ is the greatest value of $-\vecc^\T\vez/\|\vez\|_1$ over all $\vez\in\Graver(A)$ for which $\vex_k+\epsilon\vez$ (with sufficiently small $\epsilon>0$) is feasible. Rewriting this, we get
%\[
%  \Delta\geq\frac{\vecc^\T (\vex_k-\vex_{\min})}{2n\gamma}.
%\]
%Now let $\alpha\vez$ be the discrete steepest descent augmentation applied to $\vex_k$, leading to $\vex_{k+1}:=\vex_k+\alpha\vez$. Then we get
%\[
%  \vecc^\T (\vex_k-\vex_{k+1})
%  =-\alpha\vecc^\T\vez
%  =\textcolor{red}{\alpha\Delta\|\vez\|_1
%  \geq \Delta}
%  \geq \frac{\vecc^\T (\vex_k-\vex_{\min})}{2n\gamma}.
%\]
%Thus, in the LP and ILP case, we have again a factor of objective-function decrease at each augmentation that leads to a number of augmentations via Lemma \ref{Lemma: Sufficient improvement leads to polynomial algorithm} that is bounded by $O(n^2\,\gamma \,\log(\delta\,\vecc^\T (\vex_0-\vex_{\min})))$ (LP case) and by $O(n\,\gamma \,\log(\vecc^\T (\vex_0-\vex_{\min})))$ (ILP case), respectively.

%\medskip
%This completes the proof of part (a) of Theorem \ref{theorem:main}. \Halmos

\subsubsection{Proof of part (b) of of Theorem \ref{theorem:main, ILP case}: Discrete steepest descent.}

We begin with a series of lemmas. Note that the proof for the LP case follows exactly the same lines, since we do not use integrality of the components in our arguments.

\begin{lemma}\label{Lemma: Steepest descent direction is in G(A)}
  Let $\vex_k$ be a feasible solution, and let $\vez$ be an associated  steepest \-descent direction. Then there is some augmentation direction $\veg\in\Graver(A)$ from $\vex_k$, with $-\vecc^\T\veg / \|\veg\|_1 \geq -\vecc^\T\vez/\|\vez\|_1$.
\end{lemma}

\proof{Proof.}
  There is a sign-compatible representation $\vez=\sum\alpha_i\veg_i$  via elements $\veg_i\in\Graver(A)$. Observe that due to the sign-compatible representation, $\vex_k+\alpha_i\veg_i$ is also feasible for all $i$. (The components of $\vex_k+\alpha_i\veg_i$ lie between those of $\vex_k$ and $\vex_k+\vez$, implying that $\ve 0\leq\vex_k+\alpha_i\veg_i\leq\veu$.) In other words, all $\alpha_i\veg_i$ are applicable augmentations at $\vex_k$.

  It remains for us to show that there exists some index $i$ with $-\vecc^\T\veg_i/\|\veg_i\|_1\geq-\vecc^\T\vez/\|\vez\|_1$. Assume to the contrary that we have $-\vecc^\T\veg_i/\|\veg_i\|_1 < -\vecc^\T\vez/\|\vez\|_1$ for all $i$. This yields
  \begin{eqnarray*}
    -\vecc^\T\vez & = &- \sum\alpha_i\vecc^\T\veg_i\\
    & = & \sum\alpha_i\|\veg_i\|_1\frac{-\vecc^\T\veg_i}{\|\veg_i\|_1}\\[5pt]
    & < & \sum\alpha_i\|\veg_i\|_1\frac{-\vecc^\T\vez}{\|\vez\|_1}\\[5pt]
    & = & \frac{-\vecc^\T\vez}{\|\vez\|_1}\ \sum\alpha_i\|\veg_i\|_1\\[5pt]
    & = & \frac{-\vecc^\T\vez}{\|\vez\|_1}\ \|\vez\|_1\\[5pt]
    & = & -\vecc^\T\vez,
  \end{eqnarray*}
  a contradiction.
 \Halmos
\endproof

Lemma \ref{Lemma: Steepest descent direction is in G(A)} states that among all steepest-descent directions applicable at a feasible solution $\vex_k$, there is always one in $\Graver(A)$. Or in other words, the discrete steepest-descent rule is a steepest descent rule as claimed in Theorems \ref{theorem:main, ILP case} and \ref{theorem:main, LP case}, parts(c).
%Thus, we may restrict ourselves to choosing only Graver basis directions as augmenting directions.

\begin{lemma}\label{Lemma: monotonicity}
  Let $\vex_k$ be a feasible solution, let $\alpha\ves$ be a steepest-descent augmentation relative to $\vex_k$, and let $\vex_{k+1}:=\vex_k + \alpha\ves$.
  Let  $\beta\vet$ be a steepest-descent augmentation relative to $\vex_{k+1}$.
  Then we have $-\vecc^\T\ves/\|\ves\|_1\geq -\vecc^\T\vet/\|\vet\|_1$.
\end{lemma}

\proof{Proof.}
  Suppose on the contrary that $-\vecc^\T\ves/\|\ves\|_1 < -\vecc^\T\vet/\|\vet\|_1$.
  First observe that $\alpha\ves+\beta\vet$ is an applicable augmentation at $\vex_k$. Moreover, we have
  \begin{eqnarray*}
    -\vecc^\T(\alpha\ves+\beta\vet) & = & \alpha\|\ves\|_1\frac{-\vecc^\T\ves}{\|\ves\|_1}+\beta\|\vet\|_1\frac{-\vecc^\T\vet}{\|\vet\|_1}\\[5pt]
    & > & \alpha\|\ves\|_1\frac{-\vecc^\T\ves}{\|\ves\|_1}+\beta\|\vet\|_1\frac{-\vecc^\T\ves}{\|\ves\|_1}\\[5pt]
    & = & (\alpha\|\ves\|_1+\beta\|\vet\|_1)\frac{-\vecc^\T\ves}{\|\ves\|_1},\\[5pt]
    & \geq & (\|\alpha\ves+\beta\vet\|_1)\frac{-\vecc^\T\ves}{\|\ves\|_1},
  \end{eqnarray*}
  and therefore
  \[
    \frac{-\vecc^\T(\alpha\ves+\beta\vet)}{\|\alpha\ves+\beta\vet\|_1}>\frac{-\vecc^\T\ves}{\|\ves\|_1}.
  \]
  This contradicts the fact that $\ves$ was a steepest descent direction for $\vex_k$.
 \Halmos
\endproof

Lemma \ref{Lemma: monotonicity} states that the steepness of steepest-descent augmentations never  increases.

\begin{lemma}\label{Lemma: chance of sign-pattern leads to change in slope}
  Let $\vex_k$ be a feasible solution and let $\alpha_1\vez_1,\ldots,\alpha_j\vez_j$ be the following steepest-descent augmentations applied to $\vex_k$. If $\vez_1$ and $\vez_j$ do not have the same sign-pattern from $\set{\,\leq 0,\geq 0\,}^n$, then we have $-\vecc^\T\vez_1/\|\vez_1\|_1>-\vecc^\T\vez_j/\|\vez_j\|_1$.
\end{lemma}

\proof{Proof.}
  Assume on the contrary that $-\vecc^\T\vez_1/\|\vez_1\|_1 \leq -\vecc^\T\vez_j/\|\vez_j\|_1$. By monotonicity, Lemma \ref{Lemma: monotonicity}, we must have $\frac{-\vecc^\T\vez_1}{\|\vez_1\|_1}=\frac{-\vecc^\T\vez_2}{\|\vez_2\|_1}=\dots=\frac{-\vecc^\T\vez_j}{\|\vez_j\|_1}$. We conclude that
  \begin{eqnarray*}
    -\vecc^\T\textstyle\left(\sum_{i=1}^j\alpha_i\vez_i\right) & = & \textstyle\sum_{i=1}^j\alpha_i\|\vez_i\|_1\frac{-\vecc^\T\vez_i}{\|\vez_i\|_1}\\
    & = & \textstyle\sum_{i=1}^j\alpha_i\|\vez_i\|_1\frac{-\vecc^\T\vez_1}{\|\vez_1\|_1}\\
    & = & \textstyle\left(\sum_{i=1}^j\alpha_i\|\vez_i\|_1\right)\frac{-\vecc^\T\vez_1}{\|\vez_1\|_1}\\
    & > & \textstyle\left\|\sum_{i=1}^j\alpha_i\vez_i\right\|_1\frac{-\vecc^\T\vez_1}{\|\vez_1\|_1},
  \end{eqnarray*}
  since $\vez_1$ and $\vez_j$ do not have the same sign-pattern. This implies that
  \[
    \frac{-\vecc^\T\textstyle\left(\sum_{i=1}^j\alpha_i\vez_i\right)}{\|\sum_{i=1}^j\alpha_i\vez_i\|_1}>\frac{-\vecc^\T\vez_1}{\|\vez_1\|_1}.
  \]
  As by assumption $\sum_{i=1}^j\alpha_i\vez_i$ is an applicable augmentation at $\vex_k$. This contradicts the fact that $\vez_1$ was a steepest-descent augmentation for $\vex_k$.
 \Halmos
\endproof
Lemma \ref{Lemma: chance of sign-pattern leads to change in slope} states that once the sign pattern
of a steepest-descent augmentation strictly changes, the steepness must decrease.
This lemma has a surprising consequence if we only apply steepest-descent directions from $\Graver(A)$:

\begin{corollary}\label{Corollary: number of augmentations bounded by size of Graver basis}
For discrete steepest descent, no direction from $\Graver(A)$ is chosen twice as an augmenting direction. Therefore, the number of steepest-descent augmentations needed to reach an optimal solution is bounded by $|\Graver(A)|$ and thus is independent of $\veb$, $\vecc$ and the initial solution $\vez_0$.
\end{corollary}

\proof{Proof.}
  Let $\vex_k$ be a feasible solution and let $\alpha_1\vez_1,\ldots,\alpha_j\vez_j$, with $\vez_1,\ldots,\vez_j\in\Graver(A)$, be the following steepest-descent augmentations applied to $\vex_k$ with $\alpha_i$ chosen maximally in each augmentation. Moreover assume that $\vez_j=\vez_1$. Then, by Lemma \ref{Lemma: chance of sign-pattern leads to change in slope}, all intermediate augmentations must have the same sign-pattern as $\vez_1$, as otherwise the steepness of the augmentations would have dropped. As all vectors $\alpha_1\vez_1,\ldots,\alpha_j\vez_j$ have the same sign-pattern, the components of $\vex_k+\alpha_1\vez_1+\alpha_j\vez_j$ lie between the components of $\vex_k$ and $\vex_k+\sum_{i=1}^j\alpha_i\vez_i$ and therefore also between $\ve 0$ and $\veu$. Thus, $\vex_k+\alpha_1\vez_1+\alpha_k\vez_j=\vex_k+(\alpha_1+\alpha_j)\vez_1$ is a feasible solution. This contradicts the fact that $\alpha_1$ was chosen maximally.
 \Halmos
\endproof

\subsection{Proof of Lemma \ref{LP-greedy} and Theorem \ref{theorem:main, LP case}}

Let us assume that $\vex_k$ is a non-optimal feasible solution, and let $\vex_{\min}$ be an optimal solution to the LP. Then there exists a (sign-compatible) representation
\[
  \vex_{\min} -  \vex_k =\sum\alpha_i\veg_i,
\]
with $\alpha_i>0$ and with $\alpha_i\veg_i\red \vex_{\min} - \vex_k$. Moreover, due to Carath\'eodory theorem, at most $n$ summands are needed in such a representation.
Again, sign-compatibility of the representation $\vex_{\min}-\vex_k =\sum\alpha_i\veg_i$ implies that for all $i$ the vectors $\vex_{k}+\alpha_i\veg_i$ and $\vex_{\min}-\alpha_i\veg_i$ are all feasible solutions.

Let us now analyze what happens for each choice of augmentation rule:

\subsubsection{Proof of Lemma \ref{LP-greedy} for Discrete deepest descent.} We observe that
\[
  0>\vecc^\T (\vex_{\min}-\vex_k)
   =\vecc^\T \sum\alpha_i\veg_i
   = \sum\alpha_i\vecc^\T\veg_i
   \geq -n\Delta
\]
where $\Delta>0$ is the largest value of $-\alpha\vecc^\T\vez$ over all $\vez\in\Circuits(A)$ and $\alpha>0$ for which $\vex_k+\alpha\vez$ is feasible. Rewriting this, we get
\[
  \Delta\geq\frac{\vecc^\T (\vex_k-\vex_{\min})}n.
\]
Now let $\alpha\vez$ be the discrete deepest-descent augmentation applied to $\vex_k$, leading to $\vex_{k+1}:=\vex_k+\alpha\vez$. Then we get $\Delta=-\alpha\vecc^\T\vez$ and
\[
  \vecc^\T (\vex_k-\vex_{k+1})=-\alpha\vecc^\T\vez=\Delta
  \geq\frac{\vecc^\T (\vex_k-\vex_{\min})}n.
\]
Thus, we have a factor of $\beta=1/n$ of objective-function value decrease at each augmentation. Applying Lemma \ref{Lemma: Sufficient improvement leads to polynomial algorithm} with $\epsilon=1/\delta$ then yields a solution $\bar\vex$ with $|\vecc^\T(\bar\vex-\vex_{\min})|<1/\delta$ within $2n\,\log(\delta\,\vecc^\T (\vex_0-\vex_{\min})))$ many augmentations. Due to the definition of $\delta$ as the least common multiple of all subdeterminants of $A$, any vertex with an objective value of at most $\vecc^\T\bar\vex$ must be optimal. As explained right before the statement of Theorem \ref{theorem:main, ILP case}, such a vertex can be found from $\bar\vex$ in at most $n$ additional augmentations. Finally, note that once the discrete deepest-descent augmentation makes progress in objective value less than $\epsilon/(2n-2)$, we have
\[
  |\vecc^\T(\vex_k-\vex_{\min})|=\sum|\alpha_i\vecc^\T\veg_i| <(2n-2)\cdot\epsilon/(2n-2)=\epsilon.
\]
Hence, we can decide effectively when we should stop making discrete deepest-descent augmentations and should rather find a nearby vertex.

\subsubsection{Proof of part (a) of Theorem \ref{theorem:main, LP case}: Discrete Dantzig descent.}  In order to avoid zig-zagging (see example in Section 4 of \cite{Hemmecke:PSP} and the details of
how to avoid this), we must  augment to a vertex with lower objective-function value after each discrete Dantzig-descent augmentation. For this, we need at most $n$ discrete Dantzig-descent augmentations within smaller and smaller faces of $P$. Let $\vex_k$ be a vertex of the given polyhedron. Again we observe that
\[
  0>\vecc^\T (\vex_{\min}-\vex_k)
   =\vecc^\T \sum\alpha_i\veg_i
   = \sum\alpha_i\vecc^\T\veg_i
   \geq - \Delta_0 \sum\alpha_i
   \geq -n\Delta \alpha_{\max},
\]
where $\Delta_0>0$ is the greatest value of $-\vecc^\T\vez$ over all $\vez\in\Circuits(A)$ for which $\vex_k+\lambda\vez$ is feasible for sufficiently small $\lambda>0$ and where $\alpha_{\max}=\max\set{\, \alpha_i\, }$. Rewriting this, we get
\[
  \Delta_0 \geq\frac{\vecc^\T (\vex_k-\vex_{\min})}{n\alpha_{\max}}.
\]
Now, let $\alpha\vez$ be the discrete Dantzig-descent augmentation applied to $\vex_k$, leading to $\vex_{k+1}:=\vex_k+\alpha\vez$. As $\vex_k$ is a vertex and $\vez$ an edge direction, we have $\vex_k,\vez\in\frac{1}{\delta}\Z^n$. As $\vex_{k+1}$ is lies on the intersection of the half-line $\set{\,\vex=\vex_k+\lambda\vez:\lambda\geq 0\,}$ with some facet of the polyhedron, we get that $\vex_{k+1}\in\frac{1}{\delta^2}\Z^n$. Consequently, $\alpha\geq 1/\delta$.

Thus, we get
\[
  \vecc^\T (\vex_k-\vex_{k+1})=-\alpha\vecc^\T\vez=\alpha\Delta_0\geq\frac{1}{\delta}\Delta_0
  \geq\frac{\vecc^\T (\vex_k-\vex_{\min})}{n\delta \alpha_{\max}}
  \geq \frac{\vecc^\T (\vex_0-\vex_{\min})}{n\delta\gamma},
\]
where $\gamma$ is the maximum entry in any feasible solution (or, equivalently, in any vertex). Thus, we have a factor of $\beta=1/(n\delta\gamma)$ of objective-function decrease at each augmentation.
Applying Lemma \ref{Lemma: Sufficient improvement leads to polynomial algorithm} with $\epsilon=1/\delta$ then yields the desired bound on the number of augmentations $2n^2\delta\gamma\,\log(\delta\,\vecc^\T (\vex_0-\vex_{\min})))$ to reach a vertex $\bar\vex$ with $|\vecc^\T(\bar\vex-\vex_{\min})|<1/\delta$. This vertex must be optimal.

\subsubsection{Proof of part (b) of Theorem \ref{theorem:main, LP case}: Discrete steepest descent.}

The proof here follows exactly the same lines of the proof to Theorem \ref{theorem:main, ILP case}, part (c).

\subsection{Proof of Theorem \ref{Theorem: N-fold LPs and IPs are strongly poly-time solvable}.}

If we keep $A$ and $B$ fixed and let $N$ vary, the binary encoding lengths of the Graver bases of $[A,B]^{(N)}$ are bounded by a polynomial in $N$, see \cite{DeLoera+Hemmecke+Onn+Weismantel:08}. More precisely, there is a constant $g(A,B)$, the so-called \emph{Graver complexity} of $A$ and $B$ (see \cite{Hosten+Sullivant,Santos+Sturmfels,Berstein+Onn:09}), given by
\[
  \max\set{\,\|\veg\|_1:\veg\in\Graver(B\Graver(A))\,},
\]
such that
\[ \left|\Graver\left([A,B]^{(N)}\right)\right|\leq\binom{N}{g(A,B)}\left|\Graver\left([A,B]^{(g(A,B))}\right)\right|\in O\left(N^{g(A,B)}\right).
\]
This means that we can find a steepest-descent direction in $\Graver\left([A,B]^{(N)}\right)$ in time polynomial in $N$ and by Theorems \ref{theorem:main, ILP case} and \ref{theorem:main, LP case}, parts (c), the number of steepest-descent augmentations to pass from any feasible solution $\vex_0$ to an optimal solution is bounded polynomially in $N$. As the input of an $N$-fold LP/ILP contains $\Theta(N)$ integer numbers (to encode $\veb$ and $\vecc$), we can augment $\vex_0$ to optimality in strongly polynomial-time. It remains for us to show how we can find such a feasible solution $\vex_0$ in strongly polynomial-time. Note that by a shift of coordinates, we may assume without loss of generality that $\vel=\ve 0$.

To find such a feasible solution $\vex_0$, consider the extended $N$-fold LP/ILP with problem matrix
{\small
  \[
    \left(\begin{array}{cccccccccccccccc}
      B & \maO    & \maO     & I_{d_B} & -I_{d_B} & B & \maO    & \maO     & I_{d_B} & -I_{d_B} & \cdots & B & \maO    & \maO     & I_{d_B} & -I_{d_B} \\
      A & I_{d_A} & -I_{d_A} & \maO    & \maO     &   &         &          &         &          &        &   &         &          &         &          \\
        &         &          &         &          & A & I_{d_A} & -I_{d_A} & \maO    & \maO     &        &   &         &          &         &          \\
        &         &          &         &          &   &         &          &         &          & \ddots &   &         &          &         &      \\
        &         &          &         &          &   &         &          &         &          &        & A & I_{d_A} & -I_{d_A} & \maO    & \maO
    \end{array}\right).
  \]
}

\noindent This is an $N$-fold matrix composed using the matrices $\bar{A}=\begin{pmatrix}A&I_{d_A}&-I_{d_A}&\maO&\maO\end{pmatrix}$ and $\bar{B}=\begin{pmatrix}B&\maO&\maO&I_{d_B}&-I_{d_B}\end{pmatrix}$. As the right-hand side of our LP/ILP, we choose the same right-hand side vector $\veb$, and as objective vector we use a $0/1$-vector with zeros in the original components and with ones in the auxiliary components. All variables get lower bounds of $0$, and the original variables get upper bounds specified by $\veu$. Due to the special form of the matrix, we can immediately write down a feasible solution to this problem. (Simply assign zeros to the original variables and set the auxiliary slack variables according to the positive and negative parts of the right-hand side values.)

  Optimizing this special linear objective function can now be done in strongly polynomial-time, since  $\bar{A}$ and $\bar{B}$ are constant the running time of this auxiliary $N$-fold LP/ILP is bounded polynomially in $N$, but does not depend on the right-hand side or the objective vector. If the optimal value of this auxiliary LP/ILP is $0$, a feasible solution to our original problem has been found (simply drop the auxiliary components). If the optimal value is positive, our original problem is infeasible.
\Halmos
\endproof

\section{Concluding remarks} 

Through the notion of Graver basis, we have obtained extensions of classical results of Bland, Edmonds, and Karp. Our new
version applies now to the case of integer-linear programs. As a consequence we have also derived  the first-ever strongly polynomial-time 
algorithm for $N$-fold integer-linear optimization.  Our new results also show that  Kitahara-Mizuno-style bounds \cite{kitaharamizuno} hold in larger generality to include augmentations that go through the interior of the polytope and are not restricted to edges. Theorem \ref{theorem:main, ILP case} is in fact an ILP extension of those bounds too.

There are at least three interesting directions for improvement and further research. First, we remark the numbers $n$ and $2n-2$ from Carath\'eodory's and Seb\"o's theorems we used in many of the arguments can be further improved to be $n-rank(A)$ and $2(n-rank(A))-2$ respectively because the arguments can be modified to use the dimension of the kernel of $A$. Such small improvement slightly strengthens several of the results of the paper, but we leave the details to the reader.
Second, as was demonstrated in \cite{circuitdiameter,manycircuitdiameters} the estimation of number of augmentations contribute to the 
estimation of the diameters of polytopes. Third, it would be interesting to improve our strongly polynomial algorithm for $N$-fold matrices.  
Currently, the number of augmentations does not depend on $b$ and $u$, but is a polynomial of degree O({g(A,B)}) with $g(A,B)$ the 
Graver complexity (see \cite{Hosten+Sullivant,Santos+Sturmfels, Berstein+Onn:09}) Thus the degree  depends on the fixed matrices $A$ and $B$.  
It would be desirable to arrive to a lower exponent  algorithm, like the one of \cite{hemmecke-onn-romanchuk:nfold-cubic} which is not  strongly 
polynomial yet, the number of augmentations is linear in  the binary encoding of $b$ and $u$, and it is only cubic in $N$, $O(N^3)$.

% Enter the text of acknowledgments here
\section*{Acknowledgments} We are grateful for the comments from Shmuel Onn.  The first and the third author are grateful for the hospitality received by the Isaac Newton Institute for Mathematical Sciences at the University of Cambridge, UK during the program ``Polynomial Optimization''. The first author is thankful to the Technical University of Munich for the support and hospitality received during his visit in Spring 2013 as a John von Neumann professor and to the NSF for the support received through grant DMS-0914107. J. Lee was partially supported by NSF grant CMMI--1160915 and ONR grant N00014-14-1-0315.

% References here (outcomment the appropriate case)

% CASE 1: BiBTeX used to constantly update the references
%   (while the paper is being written).
%\bibliographystyle{plain} % outcomment this and next line in Case 1
\bibliographystyle{ormsv080}
\bibliography{biblioAugmentation} % if more than one, comma separated

\begin{thebibliography}{37}
\expandafter\ifx\csname natexlab\endcsname\relax\def\natexlab#1{#1}\fi
\expandafter\ifx\csname url\endcsname\relax
  \def\url#1{{\tt #1}}\fi
\expandafter\ifx\csname urlprefix\endcsname\relax\def\urlprefix{URL }\fi
\expandafter\ifx\csname urlstyle\endcsname\relax
  \expandafter\ifx\csname doi\endcsname\relax
  \def\doi#1{doi:\discretionary{}{}{}#1}\fi \else
  \expandafter\ifx\csname doi\endcsname\relax
  \def\doi{doi:\discretionary{}{}{}\begingroup \urlstyle{rm}\Url}\fi \fi

\bibitem[{Ahuja et~al.(1993)Ahuja, Magnanti, and Orlin}]{Ahuja+Magnanti+Orlin}
Ahuja, R.K., T.L. Magnanti, J.B. Orlin. 1993.
\newblock {\it Network flows\/}.
\newblock Prentice Hall Inc., Englewood Cliffs, NJ.
\newblock Theory, algorithms, and applications.

\bibitem[{B\'ar\'asz and Vempala(2010)}]{baraszvempala}
B\'ar\'asz, M., S.~Vempala. 2010.
\newblock A new approach to strongly polynomial linear programming.
\newblock http://www.cc.gatech.edu/$\sim$vempala/papers/affine.pdf.

\bibitem[{Berstein and Onn(2009)}]{Berstein+Onn:09}
Berstein, Y., S.~Onn. 2009.
\newblock The {G}raver complexity of integer programming.
\newblock {\it Annals of Combinatorics\/} {\bf 13} 289--296.

\bibitem[{Bixby et~al.(1987)Bixby, Marcotte, and {L. E. Trotter,
  Jr.}}]{BixMarTr}
Bixby, R.~E., O.~M.-C. Marcotte, {L. E. Trotter, Jr.} 1987.
\newblock Packing and covering with integral feasible flows in integral
  supply-demand networks.
\newblock {\it Math. Programming\/} {\bf 39}(3) 231--239.

\bibitem[{Bland(1974)}]{blandthesis}
Bland, R.~G. 1974.
\newblock {\it Complementary orthogonal subspaces of $n$-dimensional Euclidean
  space and orientability of matroids\/}.
\newblock Thesis (Ph.D.)--Cornell University.

\bibitem[{Bland(1976)}]{blandFlorida}
Bland, R.~G. 1976.
\newblock On the generality of network flow theory.
\newblock Presented at the ORSA/TIMS Joint National Meeting, Miami, Florida,
  Fall 1976.

\bibitem[{Bland and Jensen(1992)}]{blandjensen}
Bland, R.~G., D.~L. Jensen. 1992.
\newblock On the computational behavior of a polynomial-time network flow
  algorithm.
\newblock {\it Math. Programming\/} {\bf 54}(1, Ser. A) 1--39.

\bibitem[{Bland and Edmonds(1982)}]{blandedmonds}
Bland, R.G., J.~Edmonds. 1982.
\newblock Fast primal algorithms for totally unimodular linear programming.
\newblock Presented at XIth International Symposium on Mathematical
  Programming, Bonn, West Germany (August 1982).

\bibitem[{Bonifas et~al.(2014)Bonifas, Di~Summa, Eisenbrand, H\"ahnle, and
  Niemeier}]{determinants+diameters}
Bonifas, N., M.~Di~Summa, F.~Eisenbrand, N.~H\"ahnle, M.~Niemeier. 2014.
\newblock On sub-determinants and the diameter of polyhedra.
\newblock {\it Discrete and Computational Geometry\/} {\bf 52}(1) 102--115.

\bibitem[{Borgwardt et~al.(2014{\natexlab{a}})Borgwardt, Finhold, and
  Hemmecke}]{circuitdiameter}
Borgwardt, S., E.~Finhold, R.~Hemmecke. 2014{\natexlab{a}}.
\newblock On the circuit diameter of dual transportation polyhedra.
\newblock To appear in \emph{SIAM Discrete Math}. Available as Arxiv:1405.3184.

\bibitem[{Borgwardt et~al.(2014{\natexlab{b}})Borgwardt, J.A., and
  Finhold}]{manycircuitdiameters}
Borgwardt, S., De~Loera J.A., E.~Finhold. 2014{\natexlab{b}}.
\newblock Edges vs circuits: a hierarchy of diameters in polyhedra.
\newblock Available as ArXiv:1409.7638.

\bibitem[{Cardoso and Cl\'{\i}maco(1992)}]{cardosoclimaco}
Cardoso, D., J.~Cl\'{\i}maco. 1992.
\newblock The generalized simplex method.
\newblock {\it Operations Research Letters.\/} {\bf 12} 337--348.

\bibitem[{De~Loera et~al.(2013)De~Loera, Hemmecke, and K{\"o}ppe}]{DHKbook}
De~Loera, J.~A., R.~Hemmecke, M.~K{\"o}ppe. 2013.
\newblock {\it Algebraic and geometric ideas in the theory of discrete
  optimization\/}, {\it MOS-SIAM Series on Optimization\/}, vol.~14.
\newblock Society for Industrial and Applied Mathematics (SIAM), Philadelphia,
  PA.

\bibitem[{De~Loera et~al.(2008)De~Loera, Hemmecke, Onn, and
  Weismantel}]{DeLoera+Hemmecke+Onn+Weismantel:08}
De~Loera, J.~A., R.~Hemmecke, S.~Onn, R.~Weismantel. 2008.
\newblock {$N$}-fold integer programming.
\newblock {\it Discrete Optimization\/} {\bf 5}(2) 231--241.
\newblock In Memory of George B. Dantzig.

\bibitem[{Edmonds and Karp(1972)}]{edmondskarp}
Edmonds, J., R.~M. Karp. 1972.
\newblock Theoretical improvements in algorithmic efficiency for network flow
  problems.
\newblock {\it J. ACM\/} {\bf 19}(2) 248--264.

\bibitem[{Frank and Tardos(1987)}]{andras+eva}
Frank, A., E~Tardos. 1987.
\newblock An application of simultaneous diophantine approximation in
  combinatorial optimization.
\newblock {\it Combinatorica\/} {\bf 7}(1) 49--65.

\bibitem[{Fulkerson(1968)}]{fulk}
Fulkerson, D.~R. 1968.
\newblock Networks, frames, blocking systems.
\newblock {\it Mathematics of the {D}ecision {S}ciences, {P}art 1 ({S}eminar,
  {S}tanford, {C}alif., 1967)\/}. Amer. Math. Soc., Providence, R.I., 303--334.

\bibitem[{Graver(1975)}]{Graver:75}
Graver, J.~E. 1975.
\newblock On the foundation of linear and integer programming {I}.
\newblock {\it Mathematical Programming\/} {\bf 9} 207--226.

\bibitem[{Hemmecke(2003)}]{Hemmecke:PSP}
Hemmecke, R. 2003.
\newblock On the positive sum property and the computation of {G}raver test
  sets.
\newblock {\it Math. Programming, {S}eries {B}\/} {\bf 96} 247--269.

\bibitem[{Hemmecke et~al.(2013)Hemmecke, Onn, and
  Romanchuk}]{hemmecke-onn-romanchuk:nfold-cubic}
Hemmecke, R., S.~Onn, L.~Romanchuk. 2013.
\newblock {$N$}-fold integer programming in cubic time.
\newblock {\it Mathematical Programming\/} {\bf 137} 325--341.

\bibitem[{Hemmecke et~al.(2011)Hemmecke, Onn, and
  Weismantel}]{Hemmecke+Onn+Weismantel:oracle}
Hemmecke, R., S.~Onn, R.~Weismantel. 2011.
\newblock A polynomial oracle-time algorithm for convex integer minimization.
\newblock {\it Mathematical Programming\/} {\bf 126} 97--117.

\bibitem[{Ho{\c{s}}ten and Sullivant(2007)}]{Hosten+Sullivant}
Ho{\c{s}}ten, S., S.~Sullivant. 2007.
\newblock A finiteness theorem for {M}arkov bases of hierarchical models.
\newblock {\it J. Combin. Theory Ser. A\/} {\bf 114}(2) 311--321.

\bibitem[{Kern et~al.(2012)Kern, Faigle, and Peis}]{peis}
Kern, W., U.~Faigle, B.~Peis. 2012.
\newblock Unpublished manuscript.

\bibitem[{Kitahara and Mizuno(2013)}]{kitaharamizuno}
Kitahara, T., S.~Mizuno. 2013.
\newblock A bound for the number of different basic solutions generated by the
  simplex method.
\newblock {\it Math. Program.\/} {\bf 137}(1-2, Ser. A) 579--586.

\bibitem[{Lee(1986)}]{LeeCornell}
Lee, J. 1986.
\newblock Subspaces with well-scaled frames.
\newblock Ph.D. thesis, Cornell University.

\bibitem[{Lee(1989)}]{LeeLAA}
Lee, J. 1989.
\newblock Subspaces with well-scaled frames.
\newblock {\it Linear Algebra and its Applications\/} {\bf 114-115} 21--56.
\newblock Special Issue Dedicated to Alan J. Hoffman.

\bibitem[{Onn(2010)}]{onn:nonlinear-discrete-monograph}
Onn, S. 2010.
\newblock {\it Non-Linear Discrete Optimization\/}.
\newblock Zurich Lectures in Advanced Mathematics, European Mathematical
  Society.

\bibitem[{Orlin(1997)}]{orlin97-networkdiam}
Orlin, J.B. 1997.
\newblock A polynomial time primal network simplex algorithm for minimum cost
  flows.
\newblock {\it Math. Programming\/} {\bf 78}(2, Ser. B) 109--129.
\newblock Network optimization: algorithms and applications (San Miniato,
  1993).

\bibitem[{Rockafellar(1969)}]{rocka}
Rockafellar, R.~T. 1969.
\newblock The elementary vectors of a subspace of {$R^{N}$}.
\newblock {\it Combinatorial {M}athematics and its {A}pplications ({P}roc.
  {C}onf., {U}niv. {N}orth {C}arolina, {C}hapel {H}ill, {N}.{C}., 1967)\/}.
  Univ. North Carolina Press, Chapel Hill, N.C., 104--127.

\bibitem[{Santos and Sturmfels(2003)}]{Santos+Sturmfels}
Santos, F., B.~Sturmfels. 2003.
\newblock Higher {L}awrence configurations.
\newblock {\it J. Comb. Theory Ser. A\/} {\bf 10} 151--164.

\bibitem[{Scarf(1997)}]{Scarf:97}
Scarf, H.~E. 1997.
\newblock Test sets for integer programs.
\newblock {\it Mathematical Programming, Series A\/} {\bf 79} 355--368.

\bibitem[{Schulz and Weismantel(1999)}]{SchulzWeismantel99}
Schulz, A.S, R.~Weismantel. 1999.
\newblock An oracle--polynomial time augmentation algorithm for integer
  programming.
\newblock {\it Proceedings of the 10th Annual ACM--SIAM Symposium on Discrete
  Algorithms\/}. 967--968.

\bibitem[{Schulz and Weismantel(2002)}]{SW2002}
Schulz, A.S, R.~Weismantel. 2002.
\newblock The complexity of generic primal algorithms for solving general
  integer program.
\newblock {\it Mathematics of Operations Research\/} {\bf 27} 681 -- 692.

\bibitem[{Seb{\"o}(1990)}]{Sebo90}
Seb{\"o}, A. 1990.
\newblock Hilbert bases, {C}aratheodory's theorem and combinatorial
  optimization.
\newblock Ravi Kannan, William~R. Pulleyblank, eds., {\it IPCO\/}. University
  of Waterloo Press, 431--455.

\bibitem[{Seymour(1980)}]{seymour}
Seymour, P.~D. 1980.
\newblock Decomposition of regular matroids.
\newblock {\it J. Combin. Theory Ser. B\/} {\bf 28}(3) 305--359.

\bibitem[{Tardos(1986)}]{Tardos}
Tardos, \'E. 1986.
\newblock A strongly polynomial algorithm to solve combinatorial linear
  programs.
\newblock {\it Operations Research\/} {\bf 34}(2) pp. 250--256.

\bibitem[{Ye(2011)}]{yye}
Ye, Y. 2011.
\newblock The simplex and policy-iteration methods are strongly polynomial for
  the markov decision problem with a fixed discount rate.
\newblock {\it Math. Oper. Res.\/} {\bf 36}(4) 593--603.

\end{thebibliography}

% CASE 2: BiBTeX used to generate mypaper.bbl (to be further fine tuned)
%\input{mypaper.bbl} % outcomment this line in Case 2

\end{document}